\def\versiondate{17 May 2009}
\input math.macros
\input EPSfig.macros

\checkdefinedreferencetrue
\continuousfigurenumberingtrue
\theoremcountingtrue
\sectionnumberstrue
\forwardreferencetrue
\citationgenerationtrue
\nobracketcittrue
\hyperstrue
\initialeqmacro

\input\jobname.key
\bibsty{../../texstuff/myapalike}

\def\gh{G}  
\def\gp{\Gamma}  
\def\gpe{\gamma}  

\def\bfo{{\bf 1}}

\def\fo{{\frak F}}  
\def\wsf{{\ss WSF}}
\def\fsf{{\ss FSF}}

\def\ip#1{(\changecomma #1)}

\def\changecomma#1,{#1,\,}
\def\bigchangecomma#1,{#1,\;}
\def\leftchangecomma#1,{#1,\ }


\def\HD{{\scr H}^{(2)}}   

\def\wir{{\rm W}}   
\def\fr{{\rm F}}   
\def\sleq{\preccurlyeq}   
\def\ed#1{[#1]}  
\def\uedg(#1, #2){\ed{#1, #2}}
\def\oedg(#1){\Seq{#1}}  
\def\bp{o}   

\def\tH{H}   
\def\mto#1{\buildrel #1 \over \to}
\def\bd{\partial}
\def\im{\mathop{\rm im}}
\def\vna{{\scr R}}  
\def\HH{{\scr H}}   
\def\bops{{\scr B}}  
\def\gen#1{\langle #1 \rangle}   
\def\cbd{\delta}   
\def\supp{\mathop{\rm supp}}   
\def\wlim{\mathop{\rm wlim}}  
\def\cl{\Xi}   
\def\lE{{\ell^2(E)}}
\def\M{{\scr M}}
\def\int#1{#1^\circ}  
\def\bnd{{\ss bnd}}  
\def\topbnd{{\ss bnd}\,}  
\def\cmpl#1{#1^c}   
\def\cube{{\Bbb X}}  
\def\B{{\scr B}}   
\def\F{{\scr F}}
\def\dual#1{#1^\dagger}      
\def\pc{p_{\rm c}}
\def\pu{p_{\rm u}}
\def\iB{B}   

\def\BLPSusf{\ref b.BLPS:usf/}

\def\BLPSgip{\ref b.BLPS:gip/, hereinafter referred to as BLPS (1999)%
\def\BLPSgip{BLPS \htmllocref{\bibcode{BLPS:gip}}{(1999)}}}

\ifproofmode \relax \else\head{To appear in {\it J. Topology Anal.}}
{Version of \versiondate}\fi 
\vglue20pt

\title{Random Complexes and $\ell^2$-Betti Numbers}

\author{Russell Lyons}

\abstract{Uniform spanning trees on finite graphs and their analogues on
infinite graphs are a well-studied area. On a Cayley graph of a group, we
show that they
are related to the first $\ell^2$-Betti number of the group. Our main aim,
however, is to present the
basic elements of a higher-dimensional analogue on finite and infinite
CW-complexes, which relate to the higher $\ell^2$-Betti numbers.
One consequence is a uniform isoperimetric inequality extending work of
Lyons, Pichot, and Vassout.
We also present an enumeration similar to recent work of Duval, Klivans,
and Martin.
}

\bottomII{Primary 
 60B05, 
 20F65, 
60C05. 
Secondary
 05A15, 
 05C05. 
}
{Spanning trees,
determinant, Laplacian, CW-complex, groups, homology, harmonic, matroid.}
{Research partially supported by a Visiting Miller Research Professorship at
the University of Calif., Berkeley (2001),
NSF grants DMS-0406017 and DMS-0705518 and Microsoft Research.}

\bsection{Introduction}{s.intro}

Enumeration of spanning trees in graphs began with \ref b.Kirchhoff/.
\ref b.Cayley/ evaluated this number in the special case of a complete graph.
Cayley's theorem was extended to higher dimensions by \ref b.Kalai/, who
showed that a certain enumeration of $k$-dimensional complexes in an
$(n-1)$-dimensional simplex resulted in $n$ to the power $n - 2 \choose k$.
An extension of Kalai's result to general simplicial complexes was given by
\ref b.DKM/; an extension in a different direction was given by \ref b.Adin/.
There is more than one natural way to extend the notion of spanning tree to
higher dimensions; we choose a slightly different one than the choice of
\ref b.DKM/. The fact that both choices agree in the case of a simplex
follows from the remark on p. 341 of \ref b.Kalai/.
Our choice is more closely related to matroids and this makes such objects
exist in greater generality than those of \ref b.DKM/.
We give an enumeration result similar to that of \ref b.DKM/.

Although Kirchhoff did not state any of his results using the language of
probability, they can easily and fruitfully be stated that way.
The theory of random spanning trees, chosen uniformly from among all of
them in a given graph, began again with the papers of \ref b.Broder/ and
\ref b.Aldous:ust/.
The theory was extended to infinite graphs by \ref b.Pemantle:ust/ in
response to questions of the present author.
Since then, the theory on infinite graphs has developed significantly and
led to the discovery by \ref b.Schramm:lerw/ of the SLE processes, a major
development in contemporary probability theory.
On general infinite graphs, there are two natural and important extensions
of the uniform spanning tree measures, called free and wired uniform
spanning forest measures.
See \ref b.Lyons:bird/ for a survey and \BLPSusf\ for details.

The enumerations in higher dimensions alluded to above give different
weights to different subcomplexes, depending on the torsion of their
homology groups.
Correspondingly, the probability measures we consider are not necessarily
the uniform measures, but rather, are proportional to these same weights.
In fact, on a finite CW-complex $X$ and in a given dimension $k$, we define
two probability measures on $k$-dimensional subcomplexes of $X$;
their difference
depends on the $k$-th Betti number of $X$.
Each of these measures has free and wired extensions to infinite
CW-complexes, $X$.
This gives four measures in all.
Differences among the four depend on ($\ell^2$-)homology of $X$.
In particular, all four coincide iff the reduced $k$th $\ell^2$-homology
group of $X$ vanishes.
In case $X$ admits an action by a group $\gp$ with compact quotient
$X/\gp$, a difference among the measures can be measured by the $k$th
$\ell^2$-Betti number of $X$ with respect to $\gp$.
This leads to a uniform isoperimetric inequality.
All our measures will be determinantal, whence they satisfy various strong
properties such as negative associations (\ref b.L:det/, \ref b.BBL/).

Unfortunately, we are not able to answer analogues of some of the basic
questions answered by \ref b.Pemantle:ust/, as we lack analogues of the
algorithms of \ref b.Broder/, \ref b.Aldous:ust/, or \ref b.Wilson:gen/.

We also give a suggestive analogy to percolation theory for the case of 
dimension 1.
If it could be extended, one could resolve an important question of \ref
b.Gaboriau:invar/ relating cost to the first $\ell^2$-Betti number.
Again, an algorithm that extended one known for uniform spanning trees
would be of use; alternatively, a way to deduce topology from the
definition of a determinantal probability measure would help in this case
and that of the previous paragraph.

Most of our results were announced in \ref b.L:det/, Sec. 12.

\bsection{Determinantal Measures}{s.det}

We begin with a review of the
definitions and basic properties of determinantal probability measures that
we shall use.
In fact, we restrict ourselves to determinantal measures arising from
orthogonal projections. See \ref b.L:det/ for more details and proofs.

Let $E$ be a finite set 
and let $\B$ be a nonempty collection of subsets of $E$.
Recall that the pair $\M :=
(E, \B)$ is a {\bf matroid} with {\bf bases} $\B$ if the following exchange
property is satisfied:
$$
\forall T, T' \in \B \ \ \forall e \in T\setminus T' \ \ \exists e' \in
T'\setminus T \quad (T \setminus \{e\}) \cup \{e'\} \in \B \,.
$$
All bases have the same cardinality, called the {\bf rank} of the matroid.
In our case,
$E$ will be a set of vectors in a complex vector space and $\B$ will be the
collection of maximal linearly independent subsets of $E$, where
``maximal'' means with respect to inclusion.
Matroids of this type are called {\bf vectorial} (though in general,
one allows any field to underlie the vector space, not merely the complex
numbers).
The {\bf dual} of a matroid $\M = (E, \B)$ is the matroid $\M^\perp := (E,
\B^\perp)$, where $\B^\perp := \{E \setminus T \st T \in \B\}$. 

If $E \subset \C^s$, the usual way of representing the corresponding
vectorial matroid $\M$ is by an $(s \times E)$-matrix $M$ whose columns are
the vectors in $E$ with respect to the usual basis of $\C^s$.
One calls $M$ a {\bf coordinatization matrix} of $\M$.
In this case, the rank of the matrix $M$ equals the rank of the matroid and
a base of $\M$ is set of columns forming a basis of the column space of $M$.

For subsets $A \subseteq [1, s]$, $B \subseteq E$, let $M_{A, B}$ denote the 
matrix determined by the rows of $M$ indexed by $A$ and the columns of $M$
indexed by $B$.
Let $P_H : \lE \to \lE$ be the orthogonal projection onto the row space $H$ of
$M$.
One definition of the determinantal probability measure $\P^H$ on $\B$
corresponding to $M$ is
$$
\P^H(T) = |\det M_{A,T}|^2/\det\big(M_{A, E} (M_{A, E})^*\big)
\label e.rowform
$$
for $T \in \B$ whenever the rows indexed by $A$ form a basis of $H$,
where the superscript $*$ denotes adjoint.
(One way to see that this
defines a probability measure is to use the Cauchy-Binet formula.)
As indicated by the notation, this depends on $M$ only through $H$; this
is not hard to verify by considering a change of basis, but is immediate
from another formula, 
$$
\P^H(T) = \det [Q_H]_{T, T}
\label e.Qform
$$
for $T \in \B$, where $Q_H$ is the matrix of $P_H$ .
The representation \ref e.Qform/ has a useful extension, namely, for every
$D \subseteq E$,
$$
\P^H[D \subseteq T] = \det [Q_H]_{D, D}
\,.
\label e.fullQform
$$
In case $E$ is infinite and $H$ is a closed subspace of $\lE$, the
determinantal probability measure $\P^H$ is defined via the requirement
that \ref e.fullQform/ hold for all finite $D \subset E$.

We shall use the following theorems from \ref b.L:det/.

\procl p.card
Let $E$ be a finite set and $H$ be a subspace of $\ell^2(E)$.
Then $\P^H$ is supported on the subsets of $E$ whose cardinality equals the
dimension of $H$.
\endprocl

\procl p.dualrep
Let $E$ be a finite set.
For a subspace $H \subseteq \lE$ and its orthogonal
complement $H^\perp$, we have
$$
\all {T \in 2^E}\quad \P^{H^\perp}(E \setminus T) = \P^H(T)
\,.
$$
\endprocl

Given two probability measures $\P^1$, $\P^2$ on $2^E$, we say that {\bf
$\P^2$ stochastically dominates $\P^1$} and write $\P^1 \sleq \P^2$ if
there is a random pair $(T^1, T^2) \in 2^E \times 2^E$ with $T^i \sim
\P^i$ (meaning the law of $T^i$ is $\P^i$) and such that $T^1 \subseteq
T^2$. 
We call such a random pair a {\bf monotone coupling} of $\P^1$ and $\P^2$.
(For convenience, we are mixing the definition of stochastic
domination with a theorem of \ref b.Strassen/.)

\procl t.dominate
Let $E$ be finite or infinite and let
$H \subseteq H'$ be closed subspaces of $\lE$. Then $\P^{H} \sleq
\P^{H'}$, with equality iff $H = H'$.
\endprocl

\proof The last clause about equality was not stated in \ref b.L:det/, so
we prove it here.
If $\P^H = \P^{H'}$, then for all $e \in E$, we have $\P^H[e \in T] =
\P^{H'}[e \in T]$, i.e., $\|P_H \bfo_e\| = \|P_{H'} \bfo_e\|$.
Combining this with the assumption that $H \subseteq H'$ yields $H = H'$.
\Qed

For a set $D
\subseteq E$, recall that $\F(D)$ denotes the $\sigma$-field of events
generated by the random variable $T \cap D$. 
Define the {\bf tail} $\sigma$-field to be the
intersection of $\F(E\setminus D)$ over all finite $D$.
We say that a measure $\P$ on $2^E$ has {\bf trivial tail} if 
every event in the tail $\sigma$-field has measure either 0 or 1.

\procl t.tail
Let $E$ be infinite and let
$H$ be a closed subspace of $\lE$. The measure $\P^H$ has trivial tail.
\endprocl

\bsection{Finite CW-Complexes}{s.finite}

We consider each cell of a CW-complex $X$ to be oriented (except, of
course, the 0-cells).
Write $\cl_k X$ for the set of $k$-cells of $X$.
We identify cells with the corresponding basis elements of the chain and
cochain groups, so that $\cl_k X$ forms a basis of $C_k(X; \C)$ and $C^k(X;
\C)$.
The matrix (in this basis) of the boundary map $\bd_k = \bd_{k, X}:
C_k(X; \C) \to C_{k-1}(X; \C)$ is the matrix of incidence numbers.
In the sequel, we shall not write the coefficient group $\C$.
Recall that $Z_k(X) := \ker \bd_k$, $B_k(X) := \im \bd_{k+1}$, and $H_k(X)
:= Z_k(X)/B_k(X)$.
We also have the coboundary map
$\cbd_k = \cbd_{k, X} := \bd_{k+1}^*$ with its corresponding groups,
$Z^k(X) := \ker \cbd_k$, $B^k(X) := \im \cbd_{k-1}$, and $H^k(X)
:= Z^k(X)/B^k(X)$.

Given a finite CW-complex $X$ and a subset $T \subseteq \cl_k X$ of its
$k$-cells, write $X_T$ for the subcomplex $T \cup \bigcup_{j=0}^{k-1} \cl_j
X$. 
We call $T$ a {\bf $k$-base} if it is a base of the matroid
defined by the matrix of the boundary map
$\bd_k$, i.e., if it is maximal with $Z_k(X_T) = 0$, while we call $T$ a {\bf
$k$-cobase} if it is a base of the matroid defined by the matrix of
the coboundary map
$\cbd_k$, i.e., if it is maximal with the property that the kernel of
$\cbd_k : C^k(X_T) \to C^{k+1}(X)$ is trivial.
We remark that since $X_T$ is $k$-dimensional, $T$ is a $k$-base iff
$H_k(X_T) = 0$.

In a moment, we shall define a probability measure on the set of $k$-bases;
later, we shall define another probability measure on the set of
complements of $k$-cobases.
Before giving these probability measures, we give some examples of
$k$-bases and $k$-cobases.
If $G$ is a connected graph, then the empty set is the only
0-base, while the complement of each
vertex is a 0-cobase. The 1-bases are the spanning
trees. If $G$ and $\dual G$ form a pair of dual graphs
embedded in an orientable surface with all faces contractible, then
consider the 2-complex $X$ whose 1-skeleton is $G$ and whose 2-cells are the
faces of $G$. The 1-cobases of
$X$ are the sets $T$ of edges such that for some spanning tree $T'$
of $\dual G$, each edge in $T$ crosses an edge of $T'$ and vice versa. 
The complement of each face is a 2-base of $X$,
while the empty set is the only 2-cobase.
For another example noted by \ref b.Kalai/, let $X$ be the 5-simplex.
Its 2-bases consist of 10 triangles. Some of these 2-bases form
the usual triangulation of the projective plane using 6
vertices and 10 triangles. (This triangulation arises from the regular
icosahedron by identifying antipodal points.)

Given a set $T$ of $k$-cells and $S$ of $(k-1)$-cells, we write $\bd_{S,
T}$ for the submatrix of $\bd_k$ whose rows are indexed by $S$ and columns
by $T$.
The matrix of $\bd_k$ defines a determinantal probability measure on
the set of $k$-bases as in \ref e.rowform/:
$$
\P_k(T)
:=
\P_{k, X}(T)
:=
{\det \bd_{S, T} \bd^*_{S, T} \over \det \bd_{S, \cl_k X} \bd^*_{S, \cl_k X}}
$$
for any fixed $(k-1)$-cobase $S$.
We call this measure the {\bf $k$th lower matroidal measure on $X$}.
Also, 
if we multiply this formula by $\det \bd_{S, \cl_k X} \bd^*_{S, \cl_k X}$ and sum
over $S$, then the Cauchy-Binet formula yields 
$$
\P_k(T)
=
{\det \bd^*_{\cl_{k-1} X, T} \bd_{\cl_{k-1} X, T}
\over \sum_S \det \bd_{S, \cl_k X} \bd^*_{S, \cl_k X}}
\,.
\label e.fullcol
$$
Let $t_j(L)$ denote the order of the torsion subgroup of $H_j(X_L; \Z)$.
If we write $[G]$ for the torsion subgroup of an abelian group $G$, then in our
notation, we have $t_j(L) = |[H_j(X_L; \Z)]|$.
We now show that the measure $\P_k$ is proportional to the square of the
order of the torsion subgroup of the homology group of dimension $k-1$.
Note that if $X$ is connected and $k=1$, this shows that $\P_1$ is the
uniform measure on spanning trees since 0-dimensional homology has no
torsion; this gives a short proof of the Transfer Current Theorem of
\ref b.BurPem/.

\procl p.prtor
Let $X$ be a finite CW-complex.
For each $k$, there exists $a_k$ such that for all $k$-bases $T$ of $X$,
$$
\P_k(T)
=
a_k t_{k-1}(T)^2
\,.
$$
\endprocl

To prove this, we use a presumably well-known lemma:

\procl l.lattice
Let $V$ be a subspace of\/ $\Q^n$ of dimension $r$.
Let $B_0 \subset V \cap \Z^n$ be a set of cardinality $r$ that generates
the group $V \cap \Z^n$.
For any basis $B$ of\/ $V$ that lies in $\Z^n$, identify $B$ with the matrix
whose columns are $B$ in the standard basis of\/ $\Q^n$ and write $\gen B$ for
the subgroup of\/ $\Z^n$ generated by $B$.
Then for all such $B$, we have 
$$
\det B^* B
=
|[\Z^n/\gen B]|^2 \det B_0^* B_0
\,.
$$
\endprocl

\proof
By hypothesis on $B_0$, there exists an $r \times r$ integer
matrix $A$ such that $B = B_0 A$.
We have 
$$
\det B^* B
=
\det A^* B_0^* B_0 A
=
\det A^* \det B_0^* B_0 \det A
=
(\det A)^2 \det B_0^* B_0
\,.
$$
Also, $\Z^n/\gen{B_0} = \Z^n/(V \cap \Z^n)$ is torsion free and
$[\gen{B_0} : \gen B] = |\det A|$,
whence
$$
|[\Z^n/\gen B]|
=
|[\Z^n/\gen{B_0}]| \cdot [\gen{B_0} : \gen B]
=
[\gen{B_0} : \gen B]
=
|\det A|
\,.
$$
Comparing these identities gives the result.
\Qed

\proofof p.prtor
Chain groups have integral coefficients for the duration of this proof.
By \ref e.fullcol/, $\P_k(T)$ is proportional to 
$\det \bd^*_{\cl_{k-1} X, T} \bd_{\cl_{k-1} X, T}$.
The columns of $\bd_{\cl_{k-1} X, T}$ generate the group $B_{k-1}(X_T)$.
Thus, \ref l.lattice/ shows that
$\P_k(T)$ is proportional to 
$|[C_{k-1}(X_T)/B_{k-1}(X_T)]|^2$.
Therefore, it suffices to show that
$$
[C_{k-1}(X_T)/B_{k-1}(X_T)] = [Z_{k-1}(X_T)/B_{k-1}(X_T)]
$$
in order to complete the proof.
Let $u \in [C_{k-1}(X_T)/B_{k-1}(X_T)]$. Write $u = v +B_{k-1}(X_T)$
with $v \in C_{k-1}(X_T)$. Let $n \in \Z^+$ be such that $n u = 0$, i.e.,
$n v \in B_{k-1}(X_T)$.
Since $B_{k-1}(X_T) \subseteq Z_{k-1}(X_T)$, we have $\bd (n v) = 0$, which
implies that $\bd v = 0$, i.e., that $v \in Z_{k-1}(X_T)$.
Therefore $u \in [Z_{k-1}(X_T)/B_{k-1}(X_T)]$.
\Qed

The theorem of \ref b.Kalai/ referred to in the introduction is that 
when $X$ is an $(n-1)$-simplex and $1 \le k \le n-1$,
$$
\sum_T t_{k-1}(T)^2 = 
n^{\raise5pt\hbox{${\scriptstyle n-2 \choose \scriptstyle k}$}}
\,,
$$
where the sum is over all $k$-bases of $X$.
For example, the 2-bases in the 5-simplex mentioned earlier that correspond
to the usual triangulation of the projective plane have weight 4.
Since the projective plane can be embedded\ftnote{*}{For example,
it lies in the 4-skeleton of the 5-simplex; this skeleton is compact and
naturally embedded in the 4-sphere.} in $\R^4$, one may encounter it
when taking random 2-bases in natural 4-dimensional complexes.
We shall return to enumeration in \ref s.enum/.

From now on (except in the section on enumeration or otherwise notated),
our chain and cochain coefficients will be in $\C$.
We use the usual inner-product on $C_k(X)$, which also allows us to
identify $C_k(X)$ with $C^k(X)$.

As in \ref e.Qform/,
another form of $\P_k$ is expressed using the orthogonal projection $Q_k$
of $C_k(X)$ onto the row space of $\bd_k$, i.e., onto the space of
coboundaries $B^k(X)$. In this form, we have 
$$
\P_k(T) 
=
\det [Q_k]_{T, T}
\,.
\label e.defop
$$
Of course, $B^k(X) = Z_k(X)^\perp$.

Another natural probability measure $\tilde \P^k$ on subsets of $\cl_k X$
is given by the matrix of the coboundary map
$\cbd_k$, the determinantal
probability measure corresponding to orthogonal projection on the row space
of $\cbd_k$, i.e., the column space of $\bd_{k+1}$, which is the space of
boundaries, $B_k(X)$.
The probability measure $\P^k(T) := \P^{k, X}(T) := \tilde \P^k(\cl_k X
\setminus T)$ is the determinantal probability measure corresponding to the
subspace of $k$-cocycles, $Z^k(X) = B_k(X)^\perp$ (see \ref p.dualrep/).
We call this measure the {\bf $k$th upper matroidal measure on $X$}.
It is supported by sets of $k$-cells that are complements of $k$-cobases.
Since $B^k(X) \subseteq Z^k(X)$, it follows from \ref t.dominate/ that
the upper measure $\P^k$ stochastically dominates the lower measure $\P_k$,
with equality iff $H^k(X) = 0$.
(Note that since $X$ is finite, $H^k(X)$ is isomorphic to $H_k(X)$.)
As usual, let $b_k(X)$ denote the $k$th Betti number of $X$, the dimension
of $H_k(X)$.
By \ref p.card/, one can add $b_k(X)$ $k$-cells to a sample from $\P_k$ to
get a sample from $\P^k$.
Occasionally, we shall use the reduced Betti numbers $\tilde b_k(X)$, where
$\tilde b_k(X) = b_k(X)$ for $k > 0$, but $\tilde b_0(X) = b_0(X) - 1$ (as
long as $X \ne \emptyset$).

Recall that for a subcomplex $Y$ of $X$, one writes $C_k(X, Y) :=
C(X)/C(Y)$ and that $\bd$ is defined on the corresponding chain complex,
with kernels $Z_k(X, Y)$, images $B_k(X, Y)$, and quotients $H_k(X, Y) :=
Z_k(X, Y)/B_k(X, Y)$.
Recall also that $C^k(X, Y) = \{ u \in C^k(X) \st u \restrict C_k(Y) = 0
\}$,
that $Z^k(X, Y)$ is the kernel of $\cbd_k$ on $C^k(X, Y)$, that
$B^k(X, Y)$ is the image of $\cbd_{k-1}$ on $C^{k-1}(X, Y)$, and that
$H^k(X, Y) := Z^k(X, Y)/B^k(X, Y)$.

Thus, $T$ is the complement in $\cl_k X$ of a $k$-cobase iff $T$ is minimal
with $Z^k(X, X_T) = 0$; note that since $C^{k-1}(X, X_T) = 0$, the latter
condition is equivalent to $H^k(X, X_T) = 0$, and thus to $H_k(X, X_T) = 0$.
Because the homology sequence of the pair $(X, X_T)$ is exact, this last
condition is also equivalent to the conjunction of
the surjectivity of the natural map $H_k(X_T) \to H_k(X)$ and
the injectivity of the natural map $H_{k-1}(X_T) \to H_{k-1}(X)$.

\procl p.prtorupper
Let $X$ be a finite CW-complex.
For each $k$, there exists $a_k$ such that if $T$ is the complement of a
$k$-cobase of $X$, then
$$
\P^k(T)
=
a_k |H_k(X, X_T; \Z)|^2
\,.
$$
\endprocl

\proof
Again, for this proof, all coefficient groups not explicitly given 
are $\Z$.
An argument precisely parallel to that proving \ref p.prtor/ shows that
$\P^k(T)$ is proportional to the square of the order of the torsion
subgroup of $Z^{k+1}(X)$ modulo the image under
the map $\cbd_k$ of the $k$-cochains vanishing on $C_k(X_T)$, i.e., modulo
$B^{k+1}(X, X_T)$.
Since $X_T$ is $k$-dimensional, $C^{k+1}(X, X_T) = C^{k+1}(X)$ and 
$Z^{k+1}(X, X_T) = Z^{k+1}(X)$, whence $\P^k(T)$ is proportional to
$|[H^{k+1}(X, X_T)]|^2$.
It is well known that $|[H^{k+1}(X, X_T)]| = |[H_k(X, X_T)]|$
(e.g., see Corollary 3.3 of \ref b.Hatcher/).
Since in the present case, $H_k(X, X_T; \C) = 0$, it follows that 
$H_k(X, X_T; \Z) = [H_k(X, X_T; \Z)]$.
\Qed

Here are some simple examples.
Suppose that $X$ is the 2-complex defined by a connected graph $G$
embedded in the 2-torus, all of whose faces and edges are contractible.
Let $\dual G$ be the graph dual to $G$.
Then $\P_0$ is concentrated on the empty set, while
$\P^0$ is the law of a uniform random vertex of $G$.
The uniform spanning tree of $G$ has law $\P_1$, while the edges of $G$
that do not cross a uniform spanning tree of $\dual G$ have law $\P^1$. 
If $T \sim \P^1$, then $T$
has non-contractible cycles, but no contractible cycles.
The edges of such a $T$ generate the homology $\Z^2$ of the 2-torus.
This duality is shown in the random sample of \ref f.some-3-50-50-alt/,
where the gray edges have law $\P_1$ on a $50 \times 50$ square lattice
torus graph $G$, and those edges belonging to a cycle in $\dual G$ for
$\P^1$ are shown in black, the other edges not being shown at all.
Finally, $\P_2$ is the law of the complement of
a uniform random face of $G$ and $\P^2$ is
concentrated on the full set of all 2-cells of $X$.
We conjecture that the expected number of edges that belong to a cycle for
the law $\P^1$ on an $n \times n$ square torus graph is asymptotic to $C
n^{5/4}$ for some constant $C$; cf.\ \ref b.Kenyon/.

\efiginsnumber some-3-50-50-alt x4

In many circumstances such as the preceding paragraph, one has a pair $(X,
X^*)$ of dual cell structures on an orientable $n$-dimensional manifold;
see, e.g., Chap.~10 of \ref b.SeiThr/,
p.~84 of \ref b.RouSan/, p.~59 of \ref b.Matveev/, p.~228
of \ref b.Bryant/, or p.~25 of \ref b.Fenn/.
In this case, there are bijections $\varphi_k : \cl_k X \to \cl_{n-k} X^*$
such that
the matrix of $\cbd_{n-k, X^*}$ equals that of $\bd_{k, X}$ or its
negative.
This implies that
$\P_{k, X}$ and $\P^{n-k, X^*}$ have a
coupling $\big(T, \varphi_k[\cl_k X \setminus T]\big)$.

\bsection{Infinite CW-Complexes}{s.infinite}

When $X$ is infinite, there are natural extensions
of the probability measures $\P_k$ and $\P^k$.
We shall always assume that $X$ is locally finite unless otherwise stated.
In fact, the lower and upper measures
each have two extensions, making four measures in all.

The $k$-cells form an orthonormal basis for the Hilbert space $C^{(2)}_k(X)
:= \ell^2(\cl_k X)$,
which is identified with its dual, the space of
$\ell^2$-cochains $C_{(2)}^k(X)$.
As before, $C_k(X)$ denotes the space of $k$-chains (with complex
coefficients and finite support).
Let $Z_k(X) := \{ u \in C_k(X) \st \bd_k u = 0 \}$ and
$B_k(X) := \{ \bd_{k+1} u \st u \in C_{k+1}(X) \}$ be the usual cycle and
boundary spaces.
Let $C^k_c(X)$ denote the space of $k$-cochains that vanish off a
finite set of $k$-cells.
Let $Z^k_c(X) := \{ u \in C^k_c(X) \st \cbd_k u = 0 \}$ and
$B^k_c(X) := \{ \cbd_{k-1} u \st u \in C^{k-1}_c(X) \}$ be the cocycle and
coboundary spaces that vanish off a finite set of $k$-cells.
The measures $\P^\wir_k$, $\P^\fr_k$, $\P_\wir^k$, $\P_\fr^k$ can now be
defined as the determinantal probability measures corresponding to
orthogonal projections on, respectively, $\bar B^k_c(X)$, $Z_k(X)^\perp$,
$\bar Z^k_c(X)$, or $B_k(X)^\perp$, as in \ref e.defop/, where the bars
indicate closure in the $\ell^2$-topology:
$$\eqaln{
\P_\wir^k \leftrightsquigarrow \bar Z^k_c(X)&,\quad
\P_\fr^k \leftrightsquigarrow B_k(X)^\perp\cr
\P^\wir_k \leftrightsquigarrow \bar B^k_c(X)&,\quad
\P^\fr_k \leftrightsquigarrow Z_k(X)^\perp\cr
}$$
Those with the designation $\wir$ are called {\bf wired}, while the others
are called {\bf free}, by analogy with the
case $k=1$.
In fact, $\P^\wir_1$ is the wired (uniform) spanning forest measure,
denoted $\wsf$, while $\P^\fr_1$ is the free (uniform) spanning forest
measure, denoted $\fsf$.
For more about the terminology of free and wired, see below.
Since $B^k_c(X) \perp Z_k(X)$, we have $\bar B^k_c(X) \subseteq Z_k(X)^\perp$,
whence $\P^\wir_k \sleq \P^\fr_k$.
Since 
$Z^k_c(X) \subseteq B_k(X)^\perp$, we also have $\P_\wir^k \sleq
\P_\fr^k$.
Similarly, since $B^k_c(X) \subseteq Z^k_c(X)$, we have $\P^\wir_k \sleq
\P_\wir^k$ and since $B_k(X) \subseteq Z_k(X)$, we have $\P^\fr_k \sleq
\P_\fr^k$.
Thus, all measures stochastically dominate the wired lower measure
$\P^\wir_k$, while all are
dominated by the free upper measure $\P_\fr^k$.
Hence, all four measures 
coincide iff $\P^\wir_k = \P_\fr^k$.
We have $H_k(X) = 0$ iff $Z_k(X) = B_k(X)$ iff $Z_k(X)^\perp =
B_k(X)^\perp$ iff $\P^\fr_k = \P_\fr^k$.
Likewise, $\bar Z^k_c(X) = \bar B^k_c(X)$ iff $\P^\wir_k = \P_\wir^k$,
which is implied by (but is not equivalent to) $H^k_c(X) = 0$.

When one has a pair $(X, X^*)$ of dual cell structures on
an orientable $n$-dimensional manifold, 
$\P^\fr_{k, X}$ and $\P_\wir^{n-k, X^*}$ have a
coupling $\big(T, \varphi_k[\cl_k X \setminus T]\big)$, as do
$\P^\wir_{k, X}$ and $\P_\fr^{n-k, X^*}$.

\procl r.skel
All four $k$th matroidal measures are properly defined as long as the
$k$-skeleton of $X$ is locally finite; the $(k+1)$-skeleton of $X$ need not
be locally finite.
\endprocl

We now want to show that the free and wired measures are limits of the
kinds of measures we considered on finite complexes. 

Given a finite subcomplex $A \subset X$, write $\int A$ for the {\bf
combinatorial interior} of $A$, i.e., the set of all cells of $A$ whose
coboundary vanishes off of $A$.
Although $\int A$ is not usually a subcomplex of $X$, we shall write
$C^k(\int A)$ for the space of cochains vanishing off $\int A$.
Also, let $B^k(\int A)$ be the image of the restriction of $\cbd_{k-1}$ to
$C^{k-1}(\int A) \to C^k(A)$ and $Z^k(\int A)$ be the kernel of the restriction
of $\cbd_k$ to $C^k(\int A) \to C^{k+1}(A)$.
For the determinantal
probability measures corresponding to $B^k(\int A)$ and $Z^k(\int A)$,
write $\P_{k, \int A}$ and $\P^{k, \int A}$, respectively.
Note that both these measures give random subsets of $\cl_k A$.

\procl r.relative
Instead of working with $B^k(\int A)$ and $Z^k(\int A)$ below, we could use
instead the slightly different spaces $\bigcup Z_k(X, \overline{X \setminus
A})^\perp$ and $\bigcup B_k(X, \overline{X \setminus A})^\perp$,
respectively, where we regard elements of $C_k(X, \overline{X \setminus
A})$ as subsets of $C_k(X)$.
These may be somewhat more natural topologically, but are somewhat less
explicit and yield slightly worse inequalities.
\endprocl

We call a sequence $\Seq{A_n}$ of finite subcomplexes of $X$ an {\bf
exhaustion} if $A_n \subseteq A_{n+1}$ for each $n$ and $X = \bigcup_n
A_n$. 
For probability measures $Q_n$ on subsets of $A_n$, write
$Q = \wlim Q_n$ if for all finite $B$, the restrictions of $Q_n$ to $B$
tend to the restriction of $Q$ to $B$.

The following is straightforward to check.

\procl l.spaces
Let $X$ be a locally finite complex with an exhaustion $\Seq{A_n}$.
We have
$Z_k(X) = \bigcup_n Z_k(A_n)$,
$B_k(X) = \bigcup_n B_k(A_n)$,
$B^k_c(X) = \bigcup_n B^k(\int {A_n})$,
and
$Z^k_c(X) = \bigcup_n Z^k(\int {A_n})$,
where all four unions are increasing.
\endprocl

This gives

\procl c.limits
Let $X$ be a locally finite complex with an exhaustion $\Seq{A_n}$.
Then
$\P^\fr_{k, X} = \wlim_n \P_{k, A_n}$,
$\P_\fr^{k, X} = \wlim_n \P^{k, A_n}$,
$\P^\wir_{k, X} = \wlim_n \P_{k, \int {A_n}}$,
and
$\P_\wir^{k, X} = \wlim_n \P^{k, \int {A_n}}$.
\endprocl

\procl r.wiring
We took a direct route to limits by using
subspaces, rather than finite subcomplexes. But subcomplexes can also be
used to complete the analogy to spanning forests:
Let $A$ be a finite subcomplex of $X$.
Define a new complex $A^*$ as follows.
Let $B$ be the set of cells of $A$ whose closure intersects $\int A$ and
$\cmpl B$ the rest.
The cells of $A^*$ are those in $B$ plus
one cell $z_k$ of dimension $k$ for each $k$ with $\cl_k \cmpl B \ne
\emptyset$.
Every $k$-cell in $\cmpl B$ is identified with $z_k$ in an
orientation-preserving way.
The attaching maps among the cells of $\int A$ are the same as in $A$, but 
the others are changed.
This leads to $\cbd_k z_k = - \sum_{e \in \cl_k \int A} \cbd_k e$, which
implies that $B^k(\int A) = B^k(A^*)$ and $Z^k(\int A) = Z^k(A^*)$.
Thus, one could use $A_n^*$ in place of $\int A_n$ for the limits of \ref
c.limits/, as is done traditionally in the case $k=1$ to define the wired
uniform spanning forest.
\endprocl

For a subcomplex $A \subseteq X$, define its {\bf $k$th boundary} to be
$\bnd_k(A) := \cl_k A \setminus \cl_k \int A$.
Write $\supp u$ for the support of a chain, $u$.
Our next proposition is an analogue of the fact that all the trees in the
wired or free spanning forests of infinite connected graphs are infinite.

\procl p.fol
Suppose that $X$ is locally finite, $k \ge 1$, and $\tilde b_{k-1}(X) = 0$.
If $A$ is a finite subcomplex of $X$,
then $b_{k-1}(X_\fo \cap A) \le |\bnd_{k-1}(A)|$ a.s.\ when $\fo$ has any of
the laws $\P^\wir_k$, $\P^\fr_k$, $\P_\wir^k$, or $\P_\fr^k$.
\endprocl

\proof
Since $\P^\wir_k$ is stochastically the smallest of the four measures, it
suffices to prove the inequality for it.
Let $A$ be a finite subcomplex of $X$.
Because of the hypothesis, there is a finite subcomplex $B \subset X$ such
that every $(k-1)$-cycle of $A$ is a $(k-1)$-boundary of $B$.
In fact, we may ensure that $Z_{k-1}(A) \subseteq B_{k-1}(\int B)$ (in an
extension of our earlier notation for cochains).
By \ref c.limits/, it suffices to show for all such $B$ that when $T \sim
\P_{k, \int B}$, we have $b_{k-1}(\iB_T \cap A) \le |\bnd_{k-1}(A)|$.
Let $u \in H_{k-1}(\iB_T \cap A)$.
Since $T$ forms a basis for the vector space $B_{k-1}(\int B) = \bd_k
C_k(\int B)$, we have
$B_{k-1}(\iB_T) = B_{k-1}(\int B) \supseteq Z_{k-1}(A)$.
Therefore
$u = \bd_k w + B_{k-1}(\iB_T \cap A)$ for some $w \in C_k(\iB_T)$.
Let $y \in C_k(\iB_T \cap A)$ be the restriction of $w$ to $\iB_T \cap
A$.
Then $\supp (\bd_k w - \bd_k y) \subseteq \bnd_{k-1}(A)$ and $u = \bd_k w -
\bd_k y + B_{k-1}(\iB_T \cap A)$.
This shows that every class in $H_{k-1}(\iB_T \cap A)$ is represented
by an element of $Z_{k-1}(\bnd_{k-1}(A))$, whence $b_{k-1}(\iB_T \cap
A) = \dim H_{k-1}(\iB_T \cap A) \le \dim Z_{k-1}(\bnd_{k-1}(A))
\le \dim C_{k-1}(\bnd_{k-1}(A)) = |\bnd_{k-1}(A)|$.
\Qed

We say an infinite CW-complex $X$ has {\bf bounded degree} \msnote{Is this
a good name?} if for every $k$ the map $\bd_k$ has bounded $\ell^2$-norm.
This guarantees that the four spaces
$B_{(2)}^k(X) := \im \cbd_{k-1}$, $B^{(2)}_k(X) := \im
\bd_{k+1}$, $Z_{(2)}^k(X) := \ker \cbd_k$ and
$Z^{(2)}_k(X) := \ker \bd_k$
are well defined.
Although the first two are not necessarily closed subspaces, we do have
that $\bar B^k_c(X) = \bar B^k_{(2)}(X)$
and $\bar B_k(X) = \bar B_k^{(2)}(X)$, which is the same as $B_k(X)^\perp =
\bar B_k^{(2)}(X)^\perp$.
The corresponding statements for the kernels are not always true.
We have 
$$
C^{(2)}_k(X) = \bar B_{(2)}^k(X) \oplus Z^{(2)}_k(X) = \bar
B^{(2)}_k(X) \oplus 
Z_{(2)}^k(X)
\,,
$$
whence 
$\P_k^\wir = \P_k^\fr$ iff $\bar Z_k(X) = Z_k^{(2)}(X)$ and
$\P^k_\wir = \P^k_\fr$ iff $\bar Z^k(X) = Z^k_{(2)}(X)$.
We also deduce 
the $\ell^2$-Hodge-de Rham decomposition
$$
C^{(2)}_k(X) = \bar B_{(2)}^k(X) \oplus \bar B^{(2)}_k(X) \oplus \HD_k(X)
\,,
$$
where $\HD_k(X) := Z^{(2)}_k(X) \cap Z_{(2)}^k(X)$ is the space of harmonic
$\ell^2$-$k$-chains.
Evidently, $\HD_k(X)$ is isometrically isomorphic to $H^{(2)}_k(X) :=
Z^{(2)}_k(X)/\bar B^{(2)}_k(X)$, the reduced $k$th $\ell^2$-homology group
of $X$, which is also isometrically isomorphic to the reduced $k$th
$\ell^2$-cohomology group of $X$, $Z_{(2)}^k(X)/\bar B_{(2)}^k(X)$.
All four matroidal measures 
coincide
iff
$\HD_k(X) = 0$. In this case, we shall denote the common measure by simply
$\P_k$.

In particular, suppose that $\gp$ is a countable group acting freely on $X$
by permutation of cells and the quotient $X/\gp$ is compact. (Freeness here
means that the stabilizer of each unoriented cell consists of only the
identity of $\gp$.) In this case,
we call $X$ a {\bf cocompact $\gp$-CW-complex}.
Then $X$ has bounded degree and all the above Hilbert spaces are Hilbert
$\gp$-modules.
The $k$th $\ell^2$-Betti number of $X$ is the von Neumann dimension of
$\HD_k(X)$ with respect to $\gp$: $\beta_k(X; \gp) := \dim_\gp \HD_k(X)$.
This is 0 iff $\HD_k(X) = 0$.
For more information about $\ell^2$-homology, see \ref b.Eck/.
Note that the $\ell^2$-Betti numbers of $X$ are $\gp$-equivariant homotopy
invariants of $X$: see \ref b.CG/.

Recall that a countable group $\gp$ is {\bf amenable} if it has a {\bf
F{\o}lner exhaustion}, i.e., an increasing sequence of finite subsets $F_n$
whose union is $\gp$ such that for all finite $F \subset \gp$, we have
$\lim_{n \to\infty} |F F_n \triangle F_n|/|F_n| = 0$.
For $A \subseteq X$, write $\topbnd A$ for the {\bf topological boundary} of
$A$ in $X$.
Suppose $X$ is a $\gp$-CW-complex with finite fundamental domain $D$
and $\gp$ is amenable with F{\o}lner exhaustion $\Seq{F_n}$.
Set $A_n := F_n \bar D$. Then $\Seq{A_n}$ is an exhaustion of $X$ with
$|\cl_k \topbnd A_n|/|F_n| \to 0$ as $n \to\infty$ for each $k$.
By a theorem of \ref b.DodMat/, we
have 
$$
\lim_{n \to\infty} b_k(A_n)/|F_n| = \beta_k(X; \gp)
\label e.DM
$$
for all $k$.
\ref b.Eck:DM/ gave a simpler proof, and we shall give one that is even
further streamlined, with an extension.

Fix $k$. 
Let $\Pi_n : C_k^{(2)}(X) \to C_k^{(2)}(X)$ denote the orthogonal
projection onto $C_k(A_n)$ and $d_n(H)$ denote the ordinary trace of
$\Pi_n P_H$ for a closed subspace $H$ of $C_k^{(2)}(X)$.
\ref b.Eck:DM/ noted the following:
$$
0 \le d_n(H) \le \dim \Pi_n(H)
\,,
$$
with equality on the right if $H \subseteq C_k(A_n)$; 
$$
H = H_1 \oplus H_2 \implies d_n(H) = d_n(H_1) + d_n(H_2)
\,;
$$
and
$$
0 \le d_n(H) - |F_n| \dim_\gp H \le |\cl_k \topbnd A_n|
$$
when $H$ is $\gp$-invariant.

For example, we have that 
$$\eqaln{
\dim_\gp \bar B_c^k(X) 
&= \lim_n {d_n \big(\bar B_c^k(X)\big) \over |F_n|}
\ge \limsup_n {d_n \big(B^k(A_n)\big) \over |F_n|}
= \limsup_n {\dim B^k(A_n) \over |F_n|}
\,.
}$$
Furthermore,
$$\eqaln{
\dim_\gp Z_k(X)^\perp
&=
\lim_n {d_n(Z_k(X)^\perp) \over |F_n|}
\le
\liminf_n {\dim \Pi_n(Z_k(X)^\perp) \over |F_n|}
\le
\liminf_n {\dim \Pi_n(Z_k(A_n)^\perp) \over |F_n|}
\cr&=
\liminf_n {\dim C_k(A_n) \cap Z_k(A_n)^\perp \over |F_n|}
=
\liminf_n {\dim B^k(A_n) \over |F_n|}
\,.
}$$
On the other hand, 
$\bar B_c^k(X) \subseteq Z_k(X)^\perp$, so that we have equalities
everywhere and
$$
\dim_\gp \bar B_c^k(X) = \dim_\gp Z_k(X)^\perp
=
\lim_n {\dim B^k(A_n) \over |F_n|}
\,,
$$
which implies that
$$
\bar B_c^k(X) = Z_k(X)^\perp
\,.
\label e.lowereq
$$
An exactly parallel argument shows that
$$
\bar Z_c^k(X) = B_k(X)^\perp
\label e.uppereq
$$
with
$$
\dim_\gp B_k(X)^\perp = \lim_n {\dim Z^k(A_n) \over |F_n|}
\,.
$$
Subtracting these identities, we obtain 
$$
\beta_k(X; \gp)
=
\dim_\gp B_k(X)^\perp - \dim_\gp \bar B_c^k(X)
=
\lim_n {b_k(A_n) \over |F_n|}
\,,
$$
as desired.

Another consequence of \ref e.lowereq/ and \ref e.uppereq/ is
the following:

\procl p.DM
Suppose that $\gp$ is a countable amenable group 
and $X$ is a $\gp$-CW-complex whose $k$-skeleton is cocompact.
Then $\P^\wir_k = \P^\fr_k$ and $\P_\wir^k = \P_\fr^k$.
\endprocl

Of course, if $b_k(X) = 0$, then we also obtain that $\beta_k(X; \gp) = 0$,
a result (essentially) of \ref b.CG/.

\procl r.CG
Since $Z_k^{(2)}(X) = B^k_{(2)}(X)^\perp$, it follows that we also have
$\bar Z_k(X) = Z_k^{(2)}(X)$ in the amenable case,
whence $H_k^{(2)}(X) = \bar Z_k(X)/\bar B_k(X)$.
In the case that $X$ does not have a locally finite $k$-skeleton, \ref
b.CG/ define $\beta_k(X; \gp)$ as follows.
Consider an exhaustion of $X$ by cocompact subcomplexes $X_n$.
The inclusion of $X_m$ in $X_n$ for $m < n$ induces a homomorphism $j_{m,
n}: H_k^{(2)}(X_m) \to H_k^{(2)}(X_n)$.
Clearly $\dim_\gp \im j_{m, n}$ is decreasing in $n$, so its limit exists
and is increasing in $m$. Thus, we may define 
$$
\beta_k(X; \gp)
:=
\lim_{m \to\infty} \lim_{n \to\infty} \dim_\gp \im j_{m, n}
\,.
$$
It is easy to see that this does not depend on the exhaustion chosen.
Now in the amenable case, 
if $b_k(X) = 0$, then $Z_k(X_m) \subseteq B_k(X) = \bigcup_{n \ge m}
B_k(X_n)$, whence $\lim_{n \to\infty} \dim_\gp \im j_{m, n} = 0$, so that
$\beta_k(X; \gp) = 0$.
This is a new proof of a result of \ref b.CG/.
\endprocl


Denote the number of $k$-cells in $X/\gp$ by $f_k = f_k(X/\gp)$.
Write $\fo$ for a sample from $\P_k$.

\procl p.deg Let $\gp$ be amenable and act freely on 
a complex $X$ whose $k$-skeleton is cocompact. 
If\/ $\tilde b_{k-1}(X) = 0$, then 
the $\P_k$-expected number of $k$-cells in $\fo$ per vertex
of $X$ equals 
$$
f_{k-1}/f_0 + \sum_{j=0}^{k-2} (-1)^{k+j-1} \big(f_j - \beta_j(X; \gp)\big)/f_0 
\,.
$$
This also equals the average number of $k$-cells in $\fo$ per vertex of $X$
$\P_k$-a.s.
\endprocl

\proof
The case $k=0$ is easy, so assume that $k \ge 1$.
We use the notation above.
Let $\fo$ be a sample from the matroidal measure.
Since $X_\fo$ has no $k$-cycles,
the Euler-Poincar\'e formula yields
$$\eqaln{
\sum_{j=0}^{k-1} (-1)^j |\cl_j A_n| + &(-1)^k |\cl_k (X_\fo \cap A_n)|
\cr&=
\sum_{j=0}^k (-1)^j |\cl_j (X_\fo \cap A_n)|
\cr&=
\sum_{j=0}^k (-1)^j b_j (X_\fo \cap A_n)
\cr&=
\sum_{j=0}^{k-2} (-1)^j b_j (A_n) + (-1)^{k-1} b_{k-1} (X_\fo \cap A_n)
\,.
\label e.EP
}$$
Thus, if we divide both sides of \ref e.EP/ by $|F_n| f_0$ and use \ref
p.fol/ and \ref e.DM/, we obtain as a limit the equalities desired.
%
\Qed

Write $\cube^d$ for the natural $d$-dimensional CW-complex determined by
the tiling of $\R^d$ by a unit cube and all its translates by elements of
$\Z^d$.
The following result is suggested by duality.

\procl c.Zd The $\P_k$-probability that a given $k$-cell belongs to 
$\fo$ in $\cube^d$ is $k/d$.
\endprocl

\proof
In this case, we have $f_j = {d \choose j}$ and $\beta_j(\cube^d; \Z^d) = 0$,
whence the
$\P_k$-expected number of $k$-cells per vertex equals ${d-1 \choose k-1}$.
Since the number of $k$-cells of $\cube^d$ per vertex is ${d \choose k}$, the
result follows by symmetry, all $k$-cells having the same probability.
\Qed

We are interested in the
$\P_k$-expected number of $k$-cells per vertex of $X$ in the non-amenable
case as well.
In the case of Cayley graphs, the action of $\gp$ is not free when the edges
are undirected and there are involutions among the generators.
Since the graph case is of special interest, we give the following
result first.
For simplicity of notation, we write $\deg_\fo$ for the degree in the graph
spanned by $\fo$.

\procl p.graphs 
Let $\gh$ be the Cayley graph of a group $\gp$ with respect to a symmetric
generating set, $S$. (The edges are undirected and $S$ does not contain the
identity.)
Let $\bp$ be a vertex of $\gh$.
Let $H$ be a $\gp$-invariant closed subspace of\/ $C^{(2)}_1(\gh)$ and $\fo
\sim \P^H$.
Then 
$$
\E^H[\deg_\fo \bp] = 2\, \dim_\gp H
\,.
$$
\endprocl

\proof
Let the standard basis elements of $\ell^2(\gp \times S)$ be $\{f_{\gpe, s}
\st \gpe \in \gp, s \in S\}$.
Identify $C^{(2)}_1(\gh)$ with the range in $\ell^2(\gp \times S)$ of the
map defined by sending the edge $\uedg(\gpe, \gpe s)$ to the vector
$(f_{\gpe, s} + f_{\gpe s, s^{-1}})/\sqrt 2$.
These vectors form an orthonormal basis of the range.
Then $H$ becomes identified with a subspace $H_S$ that is not only
$\gp$-invariant, but also invariant under the involutions 
$f_{\gpe, s} \mapsto f_{\gpe s, s^{-1}}$.
Write $Q$ for the orthogonal projection of $\ell^2(\gp \times S)$ onto $H_S$.
We may choose $\bp$ to be the identity of $\gp$.
By involution invariance, we have 
$$
Q f_{\bp, s} = Q f_{s, s^{-1}}
\,.
$$
Therefore, 
$$\eqaln{
\E^H[\deg_\fo \bp] 
&=
\sum_{s \in S} \P^H\big[\uedg(\bp, s) \in \fo\big]
=
\sum_{s \in S} \|Q (f_{\bp, s} + f_{s, s^{-1}})/\sqrt 2\|^2
\cr&=
\sum_{s \in S} \| \sqrt 2 Q f_{\bp, s} \|^2
=
2 \sum_{s \in S} \ip{Q f_{\bp, s}, f_{\bp, s}}
\cr&=
2\, \dim_\gp H_S
=
2\, \dim_\gp H
\,.
\Qed
}
$$

A complex $K$ is called a {\bf $K(\gp, 1)$ CW-model} if $K$ is a
CW-complex with fundamental group equal to $\gp$ and vanishing higher
homotopy groups.
In this case, if $X$ is the universal cover of $K$,
we define
$\beta_k(\gp) := \beta_k(X; \gp)$;
it depends only on $\gp$ and not on $K$.
For example, if $k=1$ and $\gp$ is finitely presented, then $\HD_1(X)$
consists of the 1-chains that are orthogonal to both $B^1_{(2)}(X)$ and
$B_1(X)$; the latter space is the space generated by the cycles in the
Cayley graph, $G$. Hence, even when $\gp$ is not finitely presented,
$\beta_1(\gp) = \dim_\gp B^1_{(2)}(G)^\perp \cap Z_1(G)^\perp$.

\procl c.betti In any Cayley graph of a group $\gp$, we have 
$$
\E_\fsf[\deg_\fo \bp] = 2 \beta_1(\gp) + 2
\,.
$$
\endprocl

\proof
By \ref p.graphs/, we have
$$
\E_\fsf[\deg_\fo \bp] 
=
2\, \dim_\gp Z_1(G)^\perp
=
2 \beta_1(\gp) + 2\, \dim_\gp \bar B^1_{(2)}(G)
=
2 \beta_1(\gp) + 2
$$
because $\delta : C^0_{(2)}(G) \to C^1_{(2)}(G)$ is injective and $\dim_\gp
C^0_{(2)}(G) = 1$.
\Qed

This identity was extended to transitive unimodular graphs by \ref
b.LPS:msf/ (see the proof of Corollary 3.24), which depends on a definition
of \ref b.Gaboriau:HD/.

Now we extend the identity to higher dimensions.

\procl p.higher
Suppose that $\gp$ is a countable group 
and $X$ is a cocompact $\gp$-CW-complex.
Let $D$ be a fundamental domain for the action of\/ $\gp$ on $X$.
Let $H$ be a $\gp$-invariant closed subspace of\/ $C^{(2)}_k(X)$.
Then $\E^H\big[|\fo \cap D|\big] = \dim_\gp H$.
In particular, 
$\E^k_\fr\big[|\fo \cap D|\big] - \E_k^\wir\big[|\fo \cap D|\big] = 
\beta_k(X; \gp)$.
\endprocl

\proof
Let the standard basis elements of $C^{(2)}_k(X)$ be $\{f_{\gpe, e}
\st \gpe \in \gp, e \in \cl_k D\}$.
Write $Q$ for the orthogonal projection onto $H$.
Let $\bp$ be the identity of $\gp$.
Then
$$
\E^H\big[|\fo \cap D|\big]
=
\sum_{e \in \cl_k D} \P^H\big[e \in \fo\big]
=
\sum_{e \in \cl_k D} \ip{Q f_{\bp, e}, f_{\bp, e}}
=
\dim_\gp H
\,.
\Qed
$$

\procl c.model
If $K$ is a $K(\gp, 1)$ CW-model with finite $k$-skeleton
and $X$ is its universal cover with fundamental domain $D$, then on $X$, we
have $\E_k^\fr\big[|\fo \cap D|\big] - \E_k^\wir\big[|\fo \cap D|\big] =
\beta_k(\gp)$.
\endprocl

\proof
Since the higher homotopy groups of $X$ also vanish, so do its homology
groups. Thus, $\P^\fr_k = \P_\fr^k$.
By definition, $\beta_k(\gp) = \beta_k(X; \gp)$.
\Qed

We now give an extension of \ref e.DM/ to the non-amenable setting.
Our proof also gives an alternative proof that in the amenable case,
$\beta_k(X; \gp) = 0$.

\procl c.CGquant
Suppose that $\gp$ is a countable group 
and $X$ is a $\gp$-CW-complex whose $k$-skeleton is cocompact for some
fixed $k \ge 1$.
Let $D$ be a fundamental domain for the action of\/ $\gp$ on $X$.
If\/ $\tilde b_{k-1}(X) = 0$, then 
$$
\inf \left\{ {|\bnd_{k-1}(F \bar D)| \over |F|} 
\st
F \subset \gp \hbox{ is finite} \right\}
\ge
\beta_k(X; \gp)
\,.
$$
\endprocl

\proof
Let $F \subset \gp$ be finite and $A := F \bar D$.
The same reasoning that led to \ref e.EP/ shows that
$$
\sum_{j=0}^{k-1} (-1)^j |\cl_j A| + (-1)^k |\cl_k (X_\fo \cap A)|
=
\sum_{j=0}^{k-2} (-1)^j b_j (A) + (-1)^{k-1} b_{k-1} (X_\fo \cap A)
$$
when $\fo$ is a sample from any of the four matroidal measures.
Apply this to a monotone coupling
$(\fo, \fo^*)$ of $\P_\fr^k$ and $\P^\wir_k$  and subtract the resulting
equations to get
$$
|\fo \cap A| - |\fo^* \cap A| 
=
b_{k-1}(X_{\fo^*} \cap A) - b_{k-1}(X_{\fo} \cap A)
\le
|\bnd_{k-1}(A)|
\,,
$$
where we have applied \ref p.fol/ in the last step.
Therefore, 
$$
\E_\fr^k\big[|\fo \cap A|\big] - \E^\wir_{k}\big[|\fo \cap A|\big]
\le |\bnd_{k-1}(A)|
\,.
$$
The left-hand side is equal to
$$
|F| \cdot \Big(\E_\fr^k\big[|\fo \cap D|\big] - \E^\wir_{k}\big[|\fo \cap
D|\big]\Big)
=
|F| \beta_k(X; \gp)
$$
by \ref p.higher/, which gives the desired inequality.
\Qed

We immediately deduce the following inequality.

\procl c.LPV
Fix $k \ge 1$.
For a countable group $\gp$, every contractible $\gp$-CW-complex $X$
with fundamental domain $D$ and
for which $\cl_k X/\gp$ is finite satisfies
$$
\inf \left\{ {|\bnd_{k-1}(F \bar D)| \over |F|} 
\st
F \subset \gp \hbox{ is finite} \right\}
\ge
\beta_k(\gp)
\,.
$$
\endprocl

%
%
%

\ref c.CGquant/ extends Corollary 7 of \ref b.LPV/ to quasi-transitive
graphs acted on by $\gp$ and, of course, to higher dimensions.

Very interesting questions remain for the standard cubical CW-decomposition
$\cube^d$ of $\R^d$.
Recall that all four measures coincide.

\beginbullets

What is the $(k-1)$-dimensional (co)homology of the random
$k$-subcomplex? In the case $k=1$ of spanning forests, this asks how many
trees there are, the question answered by
\ref b.Pemantle:ust/.

If one takes the 1-point compactification of the random subcomplex, what is the
$k$-dimensional (co)homology? In the case of spanning forests, this asks
how many ends there are in the tree(s), the 
question answered partially by
\ref b.Pemantle:ust/ and completely by \BLPSusf.

\endbullets

Note that by translation-invariance of (co)homology and ergodicity of $\P_k$,
we have that the values of the (co)homology groups are constants a.s.

It follows trivially from the Alexander duality theorem and the results of
\ref b.Pemantle:ust/ and \BLPSusf\ that for $k=d-1$, we have $H_{k-1}(\fo)
= 0$ $\P_k$-a.s., while $\P_k$-a.s., the \v Cech-Alexander-Spanier
cohomology group $\check H^k(\fo \cup \infty)$ is $0$ for
$2 \le d \le 4$ and is (naturally isomorphic to) a direct sum of infinitely
many copies of $\Z$ for $d \ge 5$.
It also follows from the Alexander duality theorem and from equality of free
and wired limits that if $d=2k$, then the a.s.\ values of
$\check H^k(\fo \cup \infty)$ and $H_{k-1}(\fo)$ are the same (naturally
isomorphic), so that the two bulleted questions above are dual in that case.

\bsection{Analogy to Percolation}{s.aperc}

In the 1-dimensional case,
there is a suggestive analogy to phase transitions in Bernoulli 
percolation theory.
In that theory, given a connected graph $G$,
one considers for $0 < p < 1$ the random subgraph left after deletion of
each edge independently with probability $1 - p$.
A {\bf cluster} is a connected component of the remaining graph.
In the case of transitive graphs, there are two numbers $\pc, \pu \in [0,
1]$ such that if $0 < p < \pc$, then there are no infinite clusters a.s.;
if $\pc < p < \pu$, then there are infinitely many infinite clusters a.s.;
and if $\pu < p < 1$, then there is exactly 1 infinite cluster a.s. See
\ref b.HPS:merge/.

\procl p.phtr Let $\gh$ be a Cayley graph of an infinite group $\gp$
and $H$ be a $\gp$-invariant closed subspace of $C^{(2)}_1(\gh)$. 
\beginitems
\itemrm {(i)} If $H \subsetneq \bar B_{c}^1(\gh)$, then $\P^H$-a.s.\
infinitely many components of $\fo$ are finite.
\itemrm{(ii)} If $B_{c}^1(\gh) \subseteq H \subsetneq
Z_1(\gh)^\perp$, then $\P^H$-a.s.\ there are infinitely many infinite
components of $\fo$ and no finite components.
\enditems
\endprocl

\proof
Suppose that $H \subsetneq \bar B_{c}^1(\gh)$.
Since 
$
\E^H[\deg_\fo \bp] = 2\, \dim_\gp H
< 2\, \dim_\gp \bar B_{c}^1(\gh) = 2\,,
$
where $\bp \in \gp$,
it follows from Theorem 6.1 of \ref b.BLPS:gip/ that some
component is finite with positive $\P^H$-probability.
However, $\P^H$ has a trivial tail $\sigma$-field by \ref t.tail/,
which implies ergodicity of the $\gp$-action, whence this
event has probability 1. 
Now if there were only finitely many finite components, then picking a
vertex uniformly at random from their union would give a way to pick a
vertex at random in an invariant way, which is clearly impossible.
This proves part (i).

Now suppose that $B_{c}^1(\gh) \subseteq H \subsetneq
Z_1(\gh)^\perp$.
By \ref t.dominate/, we have $\P^H \sleq \fsf$. Since
$\P^H \ne \fsf$, it follows that in a monotone coupling $(\fo, \fo^*)$
of the two
measures, $A := \fo^* \setminus \fo$ is non-empty with positive
probability.
Let $e_0$ be an edge that lies in $A$ with positive probability and let $B$
be the $\gp$-orbit of $e_0$, which is necessarily infinite.
Because $\P[e \in A] = \fsf[e \in \fo] - \P^H[e \in \fo]$ and both terms on
the right-hand side are the same for all $e \in B$, we have that $\P[e \in
A]$ also is the same for all $e \in B$.
Therefore $E\big[|A|\big] = \infty$. The number of components of $\fo$ is
at least the size of $A$. Since the number of components of
$\fo$ is an invariant random variable, it is constant, whence infinite a.s.
On the other hand,
since $\wsf \sleq \P^H$, each component is infinite. This proves (ii).
\Qed

We believe that more is true, namely,
that if $H \subsetneq \bar B_{c}^1(\gh)$, then $\P^H$-a.s.\ all
components are finite.
However, there is no part (iii) in general, i.e., it is not true that for
every $\gp$-invariant $H \supsetneq Z_1(\gh)^\perp$, we have
$\P^H$-a.s.\ there is a unique infinite component, i.e., $\P^H$-a.s.\ $\fo$
is connected.
For a counter-example, let $\gp := \Z^2 * \Z^5$, a free product, and let
$\gh$ be its Cayley graph with respect to its natural generators.
We may decompose the edges of $\gh$ into those, $E_2$, that come from the
generators of $\Z^2$ and those, $E_5$, that come from $\Z^5$.
Let $H := C^{(2)}_1(E_2) + Z_1(E_5)^\perp$.
Clearly $H$ is $\gp$-invariant and strictly contains
$Z_1(\gh)^\perp$.
(One way to see the strict containment is to note that $\P^H$-a.s.\ every
edge in $E_2$ is present, while this is not true for $\fsf(\gh)$.)
However, $\P^H$ is the measure gotten by taking a sample from $\fsf(E_5)$
and adding to it all of $E_2$.
Since $\fsf(\Z^5)$ has infinitely many trees by a result of \ref
b.Pem:ust/, our claim follows.
Nevertheless, if for every $\epsilon > 0$ there were {\it some} $\gp$-invariant
$H \supset Z_1(\gh)^\perp$ with the two properties that $\dim_\gp H <
\dim_\gp Z_1(\gh)^\perp + \epsilon$ and that $\P^H$-almost every sample is
connected, then it would follow that $\beta_1(\gp)+1$ equals the cost of
$\gp$, which would answer an important question of \ref b.Gaboriau:invar/.
An analogous result {\it is} known for the free minimal spanning forest;
see \ref b.LPS:msf/.
The first property is not hard to ensure, i.e., that for every $\epsilon >
0$ there is some $\gp$-invariant $H \supset Z_1(\gh)^\perp$ with
$\dim_\gp H < \dim_\gp Z_1(\gh)^\perp + \epsilon$.
I am indebted to Vaughan Jones for the following proof of this fact.
We first prove a lemma.

\procl l.vj 
Let $\vna$ be a von Neumann algebra such that every non-0 projection in
$\vna$ has infinite rank (in the ordinary sense) and such that its
commutant $\vna'$ is a finite von Neumann algebra. Then $\vna$ has no
minimal projections.
\endprocl

\proof
Let $p \ne 0$ be a projection in $\vna$ on the Hilbert space $\HH$.
By Corollary 5.5.7 of \ref b.KadRing1/, we have $(p \vna p)' = p \vna'$.
If $p$ is minimal, then $p \vna' = \bops\big(p(\HH)\big)$ by Proposition
6.4.3 of \ref b.KadRing2/.
Let $p^\perp := I - p$.
Since $p^\perp A = A p^\perp = p^\perp A p^\perp$ for all $A \in \vna'$, it
follows that $\{A \in \vna' \st p^\perp A = 0\} = p \vna'$.
Now $\{A \in \vna' \st p^\perp A = 0\}$ is easily checked to be a two-sided
ideal in $\vna'$ that is closed in the weak operator topology.
Therefore it is equal to $q \vna'$ for some central projection $q \in
\vna'$ by Theorem 6.8.8 of \ref b.KadRing2/.
From the above, we conclude that $q \vna' = \bops\big(p(\HH)\big)$.
Since $\vna'$ is finite, it has a center-valued trace, $\tau$.
It is easily checked that $A \mapsto q \tau(A)$ is a
center-valued trace on $q \vna'$, so that $\bops\big(p(\HH)\big)$ is
finite.
This means that the rank of $p$ is finite, contradicting our assumption on
$\vna$.
\Qed

To apply this lemma, let $L(\gp)$ denote the left group von Neumann algebra
of $\gp$. By Theorem 6.7.2 of \ref b.KadRing2/, we have $L(\gp)' = R(\gp)$,
the right group von Neumann algebra.
Combining this with
Lemma 6.6.2 of \ref b.KadRing2/, we obtain $M_n\big(L(\gp)\big)' = R(\gp)
\otimes I_n$.
Every projection in $L(\gp)$ has infinite rank since $\gp$ is infinite.
Since $R(\gp)$ is finite, we deduce from \ref l.vj/ that
$M_n\big(L(\gp)\big)$ has no minimal projections.
Thus for every $\gp$-invariant closed subspace $K \subseteq
C_1^{(2)}(\gh)$, there is a $\gp$-invariant closed subspace $\{0\} \ne K'
\subsetneq K$.
Our claim follows easily from this by using $K := Z_1(\gh)^\perp$ 
and its subspaces.

\bsection{Enumeration}{s.enum}

Recall that $t_j(T) := |[H_j(X_T; \Z)]|$.
The normalizing constant $a_k$ in \ref p.prtor/ is the reciprocal of
the sum 
$$
h_{k-1}(X) 
:=
\sum t_{k-1}(T)^2
\,,
$$
where the sum is over all $k$-bases $T$ of $X$. Does this have an explicit
expression?
We answer this here, following the method of \ref b.DKM/.
Although our analogue of spanning tree is simpler than that of theirs,
their enumeration is simpler because their definition implies the
finiteness of certain homology groups.
We may clearly assume that the dimension $d$ of $X$ is equal to $k$.
We assume $d > 1$, since the case $d=1$ is the standard matrix-tree
theorem.

In this section, all coefficients of chain groups are in $\Z$ except
where otherwise indicated explicitly.
Given a set $S \subseteq \cl_k X$ of $k$-cells, let $Q_k(S)$ denote the
quotient of $Z_{k}(X)$ by $\big(Z_{k}(X) \cap B_{k}(X; \Q)\big) +
Z_{k}(X_{\cmpl S})$ and let $t_k'(S)$ denote its order,
where $\cmpl S := \cl_k X \setminus S$ denotes the set of $k$-cells that
are not in $S$.

The key lemma is: 

\procl l.key
Let $X$ be a finite CW-complex of dimension $d$, $T$ be a $d$-base of $X$,
and $S$ be a $(d-1)$-cobase of $X$. Then 
$$
|\det \bd_{S, T}| 
=
t_{d-1}(T) t_{d-2}(\cmpl S) t_{d-1}'(S) / t_{d-2}(X)
\,.
$$
\endprocl

\proof
Let $\Gamma := (X_T, X_{\cmpl S})$.
As in Proposition 4.1 of \ref b.DKM/, we have $\tH_d(\Gamma) = 0$ since
$\bd_{S, T}$ is nonsingular.
As in Proposition 4.2 of \ref b.DKM/, we also have that $|\det \bd_{S, T}|
= |\tH_{d-1}(\Gamma)|$.
The homology sequence of the pair $\Gamma$ is exact, which, since
$\tH_d(\Gamma) = 0$ and $\cl_k X_T = \cl_k X_{\cmpl S}$ for $k \le d-2$,
becomes
$$
0 \to \tH_{d-1}(X_{\cmpl S})
\mto{i_{d-1}} \tH_{d-1}(X_T)
\mto{j_{d-1}} \tH_{d-1}(\Gamma)
\mto{\bd_{d-1}} \tH_{d-2}(X_{\cmpl S})
\mto{i_{d-2}} \tH_{d-2}(X_T)
\to 0
\,.
\label e.relexact
$$
We claim that this induces an exact sequence of finite groups,
$$
0
\to \tH_{d-1}(X_T)/\ker j_{d-1}
\mto{[j_{d-1}]} \tH_{d-1}(\Gamma)
\mto{[\bd_{d-1}]} [\tH_{d-2}(X_{\cmpl S})]
\mto{[i_{d-2}]} [\tH_{d-2}(X_T)]
\to 0
\,,
\label e.torexact
$$
and that 
$$
|\tH_{d-1}(X_T)/\ker j_{d-1}|
=
t_{d-1}(T) t_{d-1}'(S)
\,.
\label e.quosize
$$
Since $H_{d-2}(X_T) = H_{d-2}(X)$, the result follows.

Any homomorphism of abelian groups restricts to an homomorphism of their
torsion subgroups; this is how we define the last two maps $[\bd_{d-1}]$
and $[i_{d-2}]$ above.
Since $\ker i_{d-2} = \im \bd_{d-1}$ is finite, it is contained in the
torsion subgroup, whence \ref e.torexact/ is exact at $[\tH_{d-2}(X_{\cmpl
S})]$.
In addition, since $\ker i_{d-2}$ contains only torsion elements, the
inverse image of
$[\tH_{d-2}(X_T)]$ also contains only torsion elements, whence $[i_{d-2}]$ is
onto.
This gives exactness of \ref e.torexact/ at $[\tH_{d-2}(X_T)]$.

Define $[j_{d-1}]: \tH_{d-1}(X_T)/\ker j_{d-1} \to \tH_{d-1}(\Gamma)$ as
the injective map induced by $j_{d-1}$. 
This gives exactness of \ref e.torexact/ at the remaining places
automatically.

It remains to prove \ref e.quosize/.
Now $\tH_{d-1}(X_{\cmpl S}) = Z_{d-1}(X_{\cmpl S})$ is free since $\dim
X_{\cmpl S} = d-1$.
We have $i_{d-1}$ is injective by exactness of
\ref e.relexact/ at $\tH_{d-1}(X_{\cmpl S})$.
Therefore $\im i_{d-1} \cap [\tH_{d-1}(X_T)] = 0$, so that we may identify
$[\tH_{d-1}(X_T)]$ with a subgroup $G$ of $K :=
\tH_{d-1}(X_T)/\im i_{d-1} = \tH_{d-1}(X_T)/\ker j_{d-1}$.
Thus, the proof will be completed once we show that $K/G$ is isomorphic to
$Q_{d-1}(S)$.
Now
$$
L := \tH_{d-1}(X_T)/ [\tH_{d-1}(X_T)]
=
Z_{d-1}(X_T) / \big(B_{d-1}(X_T; \Q) \cap Z_{d-1}(X_T) \big)
\,.
$$
Also, $Z_{d-1}(X_T) = Z_{d-1}(X)$ since
$C_{d-1}(X_T) = C_{d-1}(X)$ 
and,
since $T$ is a $d$-base, $B_{d-1}(X_T; \Q) = B_{d-1}(X; \Q)$.
Therefore, 
$$
L
=
Z_{d-1}(X) /
\big(B_{d-1}(X; \Q) \cap Z_{d-1}(X) \big)
\,.
\label e.1
$$
Since $\im i_{d-1} \cap [\tH_{d-1}(X_T)] = 0$, 
we may identify $\im i_{d-1}$ as a subgroup $M$ of $L$.
We have $L/M$ is isomorphic to $K/G$ and
$$
M
=
Z_{d-1}(X_{\cmpl S})/
\big(B_{d-1}(X; \Q) \cap Z_{d-1}(X_{\cmpl S}) \big)
\,.
\label e.2
$$
The quotient of \ref e.1/ by \ref e.2/
is isomorphic to $Q_{d-1}(S)$ because of the fact that for any group $D$
and subgroups $D_1, D_2$, we have an isomorphism between $(D/D_1)/(D_2/(D_1
\cap D_2))$ and $D/(D_1 D_2)$, where $D_2/(D_1 \cap D_2)$ is identified
with a subgroup of $D/D_1$.
%
\Qed

By straightforward applications of the Cauchy-Binet identity
as in \ref b.DKM/, we obtain the following corollary:

\procl c.count
Let $X$ be a finite CW-complex. Write 
$$
h'_k(X) 
:=
\sum t_k(\cmpl S)^2 t'_{k+1}(S)^2
\,,
$$
where the sum is over all $k$-cobases $S$ of $X$.
Then for any $(d-1)$-cobase $S$ of $X$, we have 
$$
h_{d-1}(X)
=
{t_{d-2}(X)^2 \over t_{d-2}(\cmpl S)^2 t'_{d-1}(S)^2} \det \bd_{S, X_d}
\bd^*_{S, X_d}
$$
and the product of the non-zero eigenvalues of $\bd_d \bd_d^*$
equals
$$
{h_{d-1}(X) h'_{d-2}(X) \over t_{d-2}(X)^2}
\,.
$$
\endprocl

\medbreak
\noindent {\bf Acknowledgements.}\enspace 
I am grateful to Damien Gaboriau for conversations on $\ell^2$-Betti
numbers, to Michael Larsen for discussions concerning \ref l.lattice/,
to Oded Schramm for help with the proof of \ref p.fol/, and to Vaughan Jones
for the proof of \ref l.vj/.
I also thank one of the referees for a very careful and perceptive
reading.

\def\noop#1{\relax}
\input \jobname.bbl

\filbreak
\begingroup
\eightpoint\sc
\parindent=0pt\baselineskip=10pt

Department of Mathematics,
Indiana University,
Bloomington, IN 47405-5701
\emailwww{rdlyons@indiana.edu}
{http://mypage.iu.edu/\string~rdlyons/}

\endgroup

\bye